\documentclass[12pt]{article}
\usepackage[english]{babel}
\usepackage{refcount}
\usepackage[shortlabels]{enumitem}
\usepackage{tikz}
\usepackage{amsfonts,amsmath,epsf,amssymb,amsthm,graphicx,afterpage}%,url}
\usepackage[
    draft = false,
    unicode = true,
    colorlinks = true,
    allcolors = blue,
    hyperfootnotes = true,
    citecolor = red
]{hyperref}
\usepackage{epsfig}
\graphicspath{{morefig/}{figures/}{figures00/}{figures02/}{figures05/}{figlms/}}

\DeclareGraphicsExtensions{.pdf,.png,.jpg}
\def\R{\mathbb R}
\def\K{K_5-45}

\newcommand{\tr}{\operatorname{tr}}

\DeclareMathOperator{\Cl}{Cl}
\DeclareMathOperator{\Int}{Int}

\long\def\comment#1\endcomment{}

\newtheoremstyle{mydefinition}% name
{3pt}%      Space above
{3pt}%      Space below
{\normalfont}%         Body font
{\parindent}%         Indent amount (empty = no indent, \parindent = para indent)
{\bfseries}% Thm head font
{.}%        Punctuation after thm head
{ }%     Space after thm head: " " = normal interword space;
{}%         Thm head spec (can be left empty, meaning `normal')

\theoremstyle{plain}
\newtheorem{theorem}{Theorem}[section]

\newtheorem{proposition}[theorem]{Proposition}
\newtheorem{lemma}[theorem]{Lemma}
\newtheorem{corollary}[theorem]{Corollary}

\theoremstyle{definition}
\newtheorem{remark}[theorem]{Remark}

\title{On winding numbers of $K_5$ minus an edge in the plane\footnote{I am grateful to E. Alkin, S. Lavrenchenko,  S. Melikhov, A. Miroshnikov,  E. Morozov and A. Skopenkov for useful discussions.
The research was supported by grant from Russian Science Foundation Grant № 25-21-00685, \href{https://rscf.ru/project/25-21-00685/}{https://rscf.ru/project/25-21-00685/}.
} }
\author{Garaev T. R.\footnote{Moscow Center for Continuous Mathematical Education.}  }

\begin{document}

\maketitle

\begin{abstract}
	Let $K$ be the graph on vertices $\{1, 2, 3, 4, 5\}$, having all edges except $(4, 5)$. 
	A continuous map $f\colon K\to \R^2$ is called an \emph{almost embedding} if $f$-images of non-adjacent edges are disjoint. 
	Take the winding numbers of the $f$-image of the oriented cycle $(1, 2, 3)$ around $f(4)$ and around $f(5)$. 
	We prove that the difference of these numbers equals $\pm 1$.
	This is surprising, because in other similar situations the analogous statement is wrong. 
	
\end{abstract}
\tableofcontents

\section{Main result}\label{s:mainres}

Denote by $[n]$ the set $\{1, \ldots, n\}$.
Denote by $K_n$ the complete graph with the vertex set $[n]$.
For a graph $K$ denote by $ij$ the edge joining vertices $i$ and $j$ in $K$.
We identify a graph and the body of the graph.

The main result of this paper (Theorems~\ref{t:sim-form} and~\ref{t:main}) is a two-dimensional analogue of the following elementary proposition.

\begin{proposition}[see proof of an equivalent statement in {\cite[Statement~3.6]{AMS1}}]\label{p:low-dim}
Let $g\colon K_3\to S^1$ be a continuous map such that $g(i)\notin g(jk)$ for every permutation $(i,j,k)$ of $[3]$.
Then the degree of $g$ equals $\pm 1$.
\end{proposition}

Denote by $w(\gamma, A)$ the \emph{winding number}\footnote{See definition for $\gamma$ a closed polygonal line in \cite[\S 2.3]{Sk18}, and for $\gamma$ a closed curve in 
\linebreak	
	\url{https://en.wikipedia.org/wiki/Winding_number}.}
 of a closed curve $\gamma \colon S^1\to \R^2$ around a point $A\in\mathbb{R}^2 \backslash \gamma(S^1)$.

\begin{theorem}\label{t:sim-form}
	Let $l$ be a closed curve in the plane.
	Take points $A_1,A_2,A_3$ on $l$ that split $l$ into paths $l_1,l_2,l_3$, where $A_i\notin l_i$ for each $i\in [3]$.  
	Let $O$ and $O'$ be points in the plane such that $O, O'\notin l$.
	For each $i\in [3]$ let $OA_i$ and $O'A_i$ be arcs (not necessarily line segments) joining $A_i$ to $O$ and $O'$, respectively. 
	Assume that 
	
	$\bullet$ for each $i\in [3]$ the union of arcs $OA_i$ and $O'A_i$ is disjoint from $l_i$;  
    %the part of the closed curve not containing $A_i$. 
	
	$\bullet$ for each  distinct $i, j \in [3]$ the  arcs $OA_i$ and $O'A_j$ are disjoint.
    
	Then $w(l,O)-w(l,O')= \pm1$. 
\end{theorem}

%??? remark
%The assumption in the first bullet of Theorem~\ref{t:sim-form} is analogous to the assumption $g(i)\notin g(jk)$ in Proposition~\ref{p:low-dim}.
%Both these assumptions are clearly necessary.  %???
%Condition in the second bullet in Theorem~\ref{t:sim-form} is also necessary.
%Figure~\ref{ris:adi_con} shows that the assumption in the second bullet of Theorem~\ref{t:sim-form} is also necessary.
%The assumption in the second bullet of Theorem~\ref{t:sim-form} is also necessary.

%\begin{figure}[h]
%		\centering
%		\includegraphics[scale=0.4]{adi_con.eps}
%		\caption{A counterexample to Theorem~\ref{t:sim-form} without the assumption in the second bullet: here $w(l,O)-w(l,O')=-3$.}
%	\label{ris:adi_con}
%	\end{figure}

Theorem~\ref{t:sim-form} is interesting because:

$\bullet$ it disproves \cite[Conjecture 1.6.a]{KS20};

$\bullet$ its mod~$2$ version is a special case of the van Kampen-Flores theorem, stating that $w(l,O)-w(l,O')$ is odd, while the latter theorem has no integer analogue, see Remarks \ref{r:vaKa}.\ref{r,e:vaKa-def}, \ref{r,e:vaKa-int_ana_ref}.

Below and in \S\ref{s:discuss} we place Theorem~\ref{t:sim-form} in the context of several interesting results. 
In particular, Theorem~\ref{t:sim-form} can be viewed as a theorem about the graph $K_5$ without an edge in the plane, see Theorem~\ref{t:main}.
 
%$\bullet$ a similar invariant for $K_{3,3}$ conjecturally assumes any odd value, see definition in Remark \ref{r:all_odd}.\ref{r,e:all_odd-K33}.
%\medskip
%$\bullet$ the analogous statements in higher dimensional Euclidean spaces are wrong, see \cite[Theorem 1.3]{Ni22} and Remark \ref{r:mot};

%Lemmas \ref{l:improv_exist} and \ref{l:improv_exist_last_case} correct a mistake in arXiv version 2 of this paper. 
%I found and corrected the mistake during my work on critical remarks by E. Alkin. 
%My work on his new critical remarks (sent about May, 18, 2025) is not yet completed. 
%So the arguments for Lemmas \ref{l:improv_exist} and \ref{l:improv_exist_last_case} is named `sketch of a proof'.
%The author did not find any mistakes in the proofs of Lemmas \ref{l:improv_exist} and \ref{l:improv_exist_last_case}, but they have not been verified by a referee.

%In Theorem~\ref{t:main} we give reformulation of Theorem~\ref{t:sim-form} in terms of  a restriction on certain winding numbers for almost embeddings of graph $K_5$ minus an edge into the plane.

Denote by $\K$ the graph obtained from $K_5$ by deleting the edge $45$. 
A map $f\colon \K\to\R^2$ is called an {\it almost embedding} if $f$-images of non-adjacent edges are disjoint. 

Let $K$ be a graph, and $f\colon K\to\R^2$ a continuous map.  
Denote by $j_1\ldots j_n$ the simple oriented cycle $(j_1, \ldots, j_n)$ in $K$.
Let $C$ be an oriented cycle in $K$, and $v$ a vertex of $K$ such that $f(v)\notin f(C)$. 
Set 
$$w_f(C, v):=w(f|_C, f(v)).$$

\begin{theorem}\label{t:main}
	For any continuous almost embedding $g\colon \K \rightarrow \mathbb{R}^2$ we have
    $$l(g):=w_g(123, 4)-w_g(123, 5)=\pm 1.$$
\end{theorem}

Denote by $K_{3,2}$ the complete bipartite graph with parts $[3]$ and $\{4, 5\}.$

%Theorem~\ref{t:main} is equivalent to Theorem~\ref{t:sim-form}.
Theorem~\ref{t:main} is surprising because its analogues for $K_4$ and $K_{3,3}$ are false: the corresponding invariants take every odd value, not only $\pm1$, see Remark~\ref{r:all_odd}.\ref{r,e:all_odd-cyc_K4}, \ref{r,e:all_odd-K33}.
The invariant $l(g)$ is equal to the invariant of $g|_{K_4}$ mentioned in the previous sentence and to some invariant of $g|_{K_{3,2}}$.
However this does not make Theorem~\ref{t:main} trivial, see Remark~\ref{r:apr}.\ref{r,e:apr-k5}, \ref{r,e:apr-k3}.

%Theorem~\ref{t:main} is equivalent to Theorem~\ref{t:sim-form}.

%The analogues of Theorem~\ref{t:main} for $K_4$ state that the corresponding invariant takes every odd value, not only $\pm 1$, see Remark~\ref{r:all_odd}.\ref{r,e:all_odd-cyc_K4}.
%The analogues of Theorem~\ref{t:main} for $K_{3,3}$ state that the corresponding invariant conjecturally takes every odd value, not only $\pm 1$, see Remark~\ref{r:all_odd}.\ref{r,e:all_odd-K33}.
% and $K_{3,3}$ are false, see rigorous formulation in Remarks~\ref{r:all_odd}.\ref{r,e:all_odd-cyc_K4} and~\ref{r:all_odd}.\ref{r,e:all_odd-K33} respectively. 
%A similar invariant for $K_4$ takes every odd value, and a similar invariant for $K_{3,3}$ is conjectured to take every odd value.
%See definitions in Remark~\ref{r:all_odd}.\ref{r,e:all_odd-cyc_K4} and~\ref{r:all_odd}.\ref{r,e:all_odd-K33}; 

The motivation, discussion and proof of Theorem~\ref{t:main} assuming Lemma \ref{p:GreatLemma2} are given in \S\ref{s:discuss}.
The proof of Lemma~\ref{p:GreatLemma2} assuming lemmas from \S\ref{s:lemmas} is given in \S\ref{s:proof}. 

\section{Discussion of Theorem \ref{t:sim-form}}  \label{s:discuss}

None of the remarks in this text are formally used in the proof of Theorem~\ref{t:sim-form}.
%See a comment Remark \ref{r:almost_emb}.\ref{r,e:almost_emb-PL-Con}

A map $g\colon K \rightarrow \R^2$ of a graph $K$ is said to be \textbf{piecewise linear} if there is a  subdivision $K'$ of $K$ such that the corresponding map $g'\colon K'\to \R^2$ is linear on any edge of $K'$.
We write 'PL' instead of 'piecewise linear'.

We define an almost embedding for any graph.
For a graph $\K$ this definition is different from but equivalent to the one given in \S 1.
A map $f\colon K\to \R^2$ of a graph $K$ is called an \textbf{almost embedding} if $f(\alpha)\cap f(\beta)=\varnothing$ for any two non-adjacent simplices (i.e. vertices or edges) $\alpha, \beta \subset K$.
In other words, if

$\bullet$ the images of non-adjacent edges are disjoint,

$\bullet$ the image of a vertex is not contained in the image of any edge non-adjacent to this vertex, 

$\bullet$ the images of distinct vertices are distinct.

We write 'embedding' and 'almost embedding' instead of 'PL or continuous embedding' and 'PL or continuous almost embedding' respectively, because of Remark \ref{r:almost_emb}.\ref{r,e:almost_emb-PL-Con}.

%In this text for an almost embedding $g\colon \K\to\R^2$ denote $l(g):=w_g(123, 4)-w_g(123, 5)$.

\begin{remark}[almost embeddings]\label{r:almost_emb}	
	\begin{enumerate}[(a)]
		\item\label{r,e:almost_emb-other_branches} Almost embeddings naturally appear in topological graph theory, in combinatorial geometry, in topological combinatorics, and in studies of embeddings (of graphs in surfaces, and of hypergraphs in higher-dimensional Euclidean space).
		See more motivations in \cite[\S 1, ‘Motivation and background’]{ST17}, \cite[\S6.10, ‘Almost embeddings, $\mathbb{Z}_2$- and $\mathbb{Z}$-embeddings’]{Sk}.
		
		\item\label{r,e:almost_emb-PL-Con}
		The property of being an almost embedding is stable, i.e., is preserved under small enough perturbation of a map (as opposed to the property of being an embedding). 
		Thus any continuous almost embedding can be approximated by a PL almost embedding. 
		For this reason Theorem \ref{t:main} is equivalent to the same statement for PL almost embeddings. 
		However even for this case Theorem \ref{t:main} is not obvious. 
		
		\item\label{r,e:almost_emb-his_alg} 
		An algebraic version of almost embeddings ($\mathbb{Z}_2$-embeddings) appeared in the 1930s and has been actively studied in graph theory since the 2000s. 
		See e.g. surveys \cite{SS13}, \cite[\S6.10 `Almost embeddings, $\mathbb{Z}_2$- and $\mathbb{Z}$-embeddings']{Sk}, and the papers \cite{FK19}, \cite{Ky16}.
		The invariant $l(g)$ and the similar invariants for graphs $K_4$ and $K_{3,3}$ defined in Remarks \ref{r:all_odd}.\ref{r,e:all_odd-cyc_K4},\ref{r,e:all_odd-K33} assume only odd values for $\mathbb{Z}_2$-embeddings $g$. 
		These are proved similarly to the analogues of these results for almost embeddings, proved in Remarks \ref{r:vaKa}.\ref{r,e:vaKa-def},\ref{r,e:vaKa-kam_k33}, \ref{r:all_odd}.\ref{r,e:all_odd-cyc_K4}.
		Presumably the analogue of Theorem~\ref{t:main}  for $\mathbb{Z}_2$-embeddings is incorrect.
		We conjecture that the analog of Theorem~\ref{t:main}  for $\mathbb{Z}$-embeddings (defined in \cite[\S 1.1]{Sk21}) is correct.
		
		\item\label{r,e:almost_emb-emb} The analogue of Theorem~\ref{t:main} for embeddings instead of almost embeddings is much simpler (and is close to the Jordan Curve Theorem). 
	\end{enumerate}
\end{remark}

\begin{remark}[motivation]\label{r:mot}
A \emph{hypergraph} is a higher-dimensional analog of a graph: together with edges joining pairs of points one considers triangles spanned by triples of points, etc., see definition e.g. in \cite[\S 6.3]{Sk}.   
A classical problem in topology, combinatorics and computer science is to find criteria for realizability (and algorithms recognizing realizability) of hypergraphs in Euclidean space of given dimension $d$. 

Such a criterion was obtained in 1930s-1960s by classical figures in topology, see survey \cite[\S4, \S5]{Sk06}.
The criterion is stated in terms of certain configuration space, yields many specific corollaries, and works for $2d\ge3k+3$, where $k$ is the dimension of the hypergraph, see survey \cite[\S5]{Sk06}. 
A polynomial algorithm based on this criterion was obtained in 2013 \cite{CKV}. 
The non-existence of a polynomial algorithm for $2d<3k+2$ was announced in 2019 by 
Marek Filakovsk\'y, Ulrich Wagner and Stephan Zhechev \cite{FWZ}. 
A mistake was found in 2020 by Arkadiy Skopenkov (and recognized by the authors). 
The mistake was that in a higher-dimensional analog of Theorem \ref{t:main} for embeddings certain linking number can assume value distinct from $\pm1$. 
In 2020 Roman Karasev and Arkadiy Skopenkov showed that this linking number for almost embeddings assumes any odd value \cite[Theorem 1.5]{KS20}.

The analogous (to the linking number) invariants for graphs in the plane are $l(g)$, and the similar invariants for the graphs $K_4$ and $K_{3,3}$ defined in Remarks \ref{r:all_odd}.\ref{r,e:all_odd-cyc_K4},\ref{r,e:all_odd-K33}.  
%We conjecture that these invariants assume any odd values, see Remarks \ref{r:all_odd}.\ref{r,e:all_odd-cyc_K4}, \ref{r,e:all_odd-K33}.
These analogues and other similar invariants assume \emph{any} odd values, see the first two bullets after Theorem~\ref{t:main}, and 
Remarks \ref{r:oth_inv}.\ref{r,e:oth_inv-cyc_wu},\ref{r,e:oth_inv-tri_wu}. %(or assume any integer value in Remark \ref{r:all_odd}.\ref{r,e:all_odd-diff} or conjecturally any even value in Remark \ref{r:oth_inv}.\ref{r,e:oth_inv-wu_diff}).
The invariant $l(g)$ and its analogues assume {\it only} odd values, see proofs in Remarks \ref{r:all_odd}.\ref{r,e:all_odd-cyc_K4},  \ref{r:vaKa}.\ref{r,e:vaKa-def},\ref{r,e:vaKa-kam_k33}.

%Theorem \ref{t:main} shows that the analogs of \cite[Theorem 1.5]{KS20} and conjectures from Remarks \ref{r:all_odd}.\ref{r,e:all_odd-cyc_K4},\ref{r,e:all_odd-K33} for $l(g)$ are false. 

\end{remark}

A homotopy $F\colon K\times [0, 1]\to\R^2$ is called an {\it (almost) isotopy}, if for any $t\in [0, 1]$ the map $F|_{K\times t}$ is an (almost) embedding.

\begin{remark}[some invariants of (almost) embeddings]\label{r:all_odd}
In this remark we introduce some invariants of an (almost) embedding up to (almost) isotopy.

	\begin{enumerate}[(a)]	
	\item\label{r,e:all_odd-diff}  
	Take an oriented cycle $C$ in a graph $K$, and some vertex $v$ in $K\backslash C$.
	For an (almost) embedding $g\colon K\to \R^2$ the integer $w_g(C, v)$ is an (almost) isotopy invariant of $g$. 
	
	For some $K$ there is an integer that is the value of this invariant for some almost embedding, but not for any embedding.
	E.g., for graph $123\sqcup \{4\}$ the number $w_g(123, 4)$ for an almost embedding $g\colon 123\sqcup \{4\}\to \R^2$ can be any integer, but for any embedding $g\colon 123\sqcup \{4\}\to \R^2$ we have $w_g(123, 4)\in \{-1, 0, 1\}$ (the latter is close to Jordan Curve Theorem). 
	Analogous statement holds if we replace the graph $123\sqcup \{4\}$ by the graph $K_4$. 
	
	For any integer $k$ there is an almost embedding $g\colon \K \to \R^2$ such that $w_g(123, 4)=k$, see Figure \ref{ris:winding}. 
	However for any embedding $g\colon \K \to \R^2$ we have $w_g(123, 4)\in \{-1, 0, 1\}$.
	
	\begin{figure}[h]
		\centering
		\includegraphics[scale=0.8]{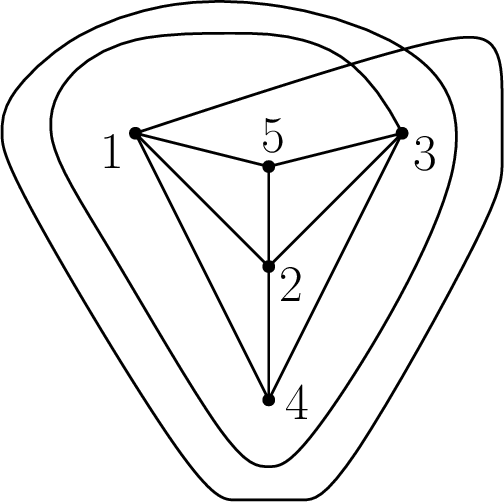}
		\caption{An almost embedding $g\colon \K \to \R^2$ such that $w_g(123, 4)=3$.}
		\label{ris:winding}
	\end{figure}
	
	\item\label{r,e:all_odd-cyc_K4}
	For $j\in [4]$ denote by $C_j$ the oriented cycle in $K_4$ obtained by deleting $j$ from $1234$.  
	Radon theorem for the plane \cite[Lemma 2.2.3]{Sk18} implies that for any almost embedding $g\colon K_4\to \R^2$ the integer $$l(g):=w_g(C_4, 4)+w_g(C_3, 3)+w_g(C_2, 2)+w_g(C_1, 1)$$ is odd, see survey \cite[Theorem 5.2]{ABM+}. 
	For any integers $n_1,n_2,n_3,n_4$ whose sum is odd there is an almost embedding $f\colon K_4\to\R^2$ such that $w_f(C_j, j)=n_j$ for every $j=1,2,3,4$ \cite[Theorem 2]{AM25}.
		
	\item\label{r,e:all_odd-K33} Take an edge $ab$ of $K_{3, 3}$. %turgor 
	Denote by $C$ somehow oriented cycle $K_{3,3}-a-b$ of length $4$. 
	It is well-known that for any almost embedding $g\colon K_{3,3}-ab \rightarrow \mathbb{R}^2$ the integer $l(g):=w_g(C, a)-w_g(C, b)$ is odd, see Remark \ref{r:vaKa}.\ref{r,e:vaKa-kam_k33}.
	We conjecture that for any integer $k$ there is an almost embedding $g\colon  K_{3, 3}- ab\to \mathbb{R}^2$ such that $l(g)=2k+1$.
	The idea was proposed by A. Lazarev, but the proof has not been published at the time of writing. 
	\end{enumerate}
\end{remark}

Denote by $K_{3,2}\subset \K$ the complete bipartite graph with parts $[3]$ and $\{4, 5\}.$
For $j\in \{4, 5\}$ denote by $[3]*j$ the complete bipartite graph with parts $[3]$ and $\{j\}$.

\begin{remark}[other invariants of (almost) embeddings]\label{r:oth_inv}
	\begin{enumerate}[(a)]
	\item\label{r,e:oth_inv-cyc_wu}
	%Denote by $l_m$ for $m\in [3]$ the edge joining two of the vertices distinct from $m$.
	The  {\it cyclic Wu number} of a map $g\colon K_3\to \mathbb{R}^2$ is defined to be the number of revolutions in the following rotation of vector:  

	from $\overrightarrow{g(1)g(2)}$ to $\overrightarrow{g(1)g(3)}$, as the second point of the vector moves along $g|_{23}$, then

	from $\overrightarrow{g(1)g(3)}$ to $\overrightarrow{g(2)g(3)}$, as the first point of the vector moves along $g|_{12}$, then

	from $\overrightarrow{g(2)g(3)}$ to $\overrightarrow{g(2)g(1)}$, as the second point of the vector moves along $g|_{31}$, then

	from $\overrightarrow{g(2)g(1)}$ to $\overrightarrow{g(3)g(1)}$, as the first point of the vector moves along $g|_{23}$, then

	from $\overrightarrow{g(3)g(1)}$ to $\overrightarrow{g(3)g(2)}$, as the second point of the vector moves along $g|_{12}$, then

	from $\overrightarrow{g(3)g(2)}$ to $\overrightarrow{g(1)g(2)}$, as the first point of the vector moves along $g|_{31}$. 

	This equals twice the (non-integer) number of revolutions in the first three rotations above.
	The cyclic Wu number is an (almost) isotopy invariant of an (almost) embedding $g\colon K_{3}\to \mathbb{R}^2$. 

	\begin{figure}[h]
	\centering
	\includegraphics[scale=0.1]{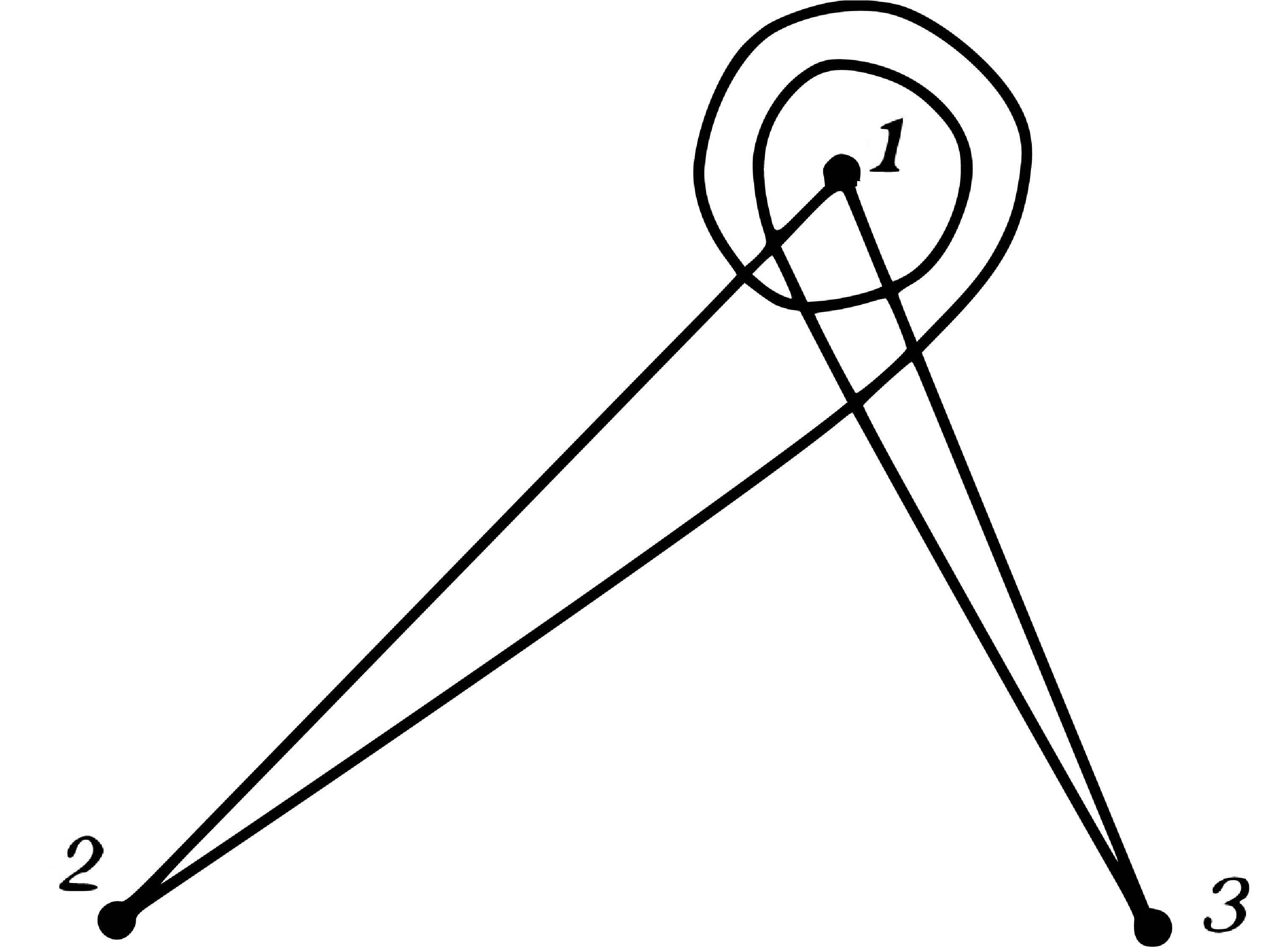}
	\caption{An almost embedding $g\colon K_{3}\to \mathbb{R}^2$ whose cyclic Wu number equals 5}
	\label{ris:K_3}
	\end{figure}

	The cyclic Wu number is odd for any almost embedding. 
	For any integer $k$ there is an almost embedding $g\colon K_{3}\to \R^2$ whose cyclic Wu number equals  $2k+1$, see Figure \ref{ris:K_3}.

	For any embedding $g\colon K_{3}\to \R^2$ the cyclic Wu number equals $\pm 1$ (this is close to Jordan Curve Theorem).

	%The cyclic Wu number is similar to, but distinct from the \emph{degree} of a closed curve. 

	\item\label{r,e:oth_inv-tri_wu}
	Denote $K_{3, 1} := [3]*4$. 
	The {\it triodic  Wu number} $\tr(g)$ of a map $g\colon K_{3,1} \to \mathbb{R}^2$ is defined to be twice the number of revolutions in the following rotation of vector: 

	from $\overrightarrow{g(1)g(2)}$ to $\overrightarrow{g(1)g(3)}$, as the second point of the vector moves along $g|_{24\cup 43}$, then

	from $\overrightarrow{g(1)g(3)}$ to $\overrightarrow{g(2)g(3)}$, as the first point of the vector moves along $g|_{14\cup 42}$, then

	from $\overrightarrow{g(2)g(3)}$ to $\overrightarrow{g(2)g(1)}$, as the second point of the vector moves along $g|_{34\cup 41}$.

	The triodic Wu number is an (almost) isotopy invariant of an (almost) embedding $g\colon K_{3, 1}\to \mathbb{R}^2$. 

	\begin{figure}[h]
		\centering
		\includegraphics[scale=0.8]{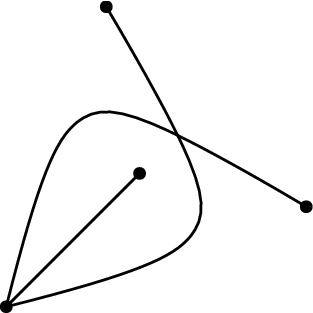}
		\caption{An almost embedding $g\colon K_{3, 1}\to \mathbb{R}^2$ whose triodic Wu number equals 3}
		\label{ris:K_3,1}
	\end{figure}  

	The triodic Wu number is odd for any almost embedding. 
	For any integer $k$ there is an almost embedding $g\colon K_{3, 1}\to \R^2$ whose triodic Wu number equals $2k+1$, see Figure \ref{ris:K_3,1}. 

	For any PL embedding, and apparently for any continuous embedding $g\colon K_{3, 1}\to \R^2$ the triodic Wu number equals $\pm 1$ (for PL embedding this is proved by induction on the number of segments in $g(ij)$ for $i\in [3]$). 
 
	\end{enumerate}
\end{remark}

\begin{remark}[on non-triviality of Theorem~\ref{t:main}]\label{r:apr}
	Denote by $g\colon \K\to\R^2$ an almost embedding.
	
	\begin{enumerate}[(a)]
	\begin{figure}[h]
	\centering
	\includegraphics[scale=1]{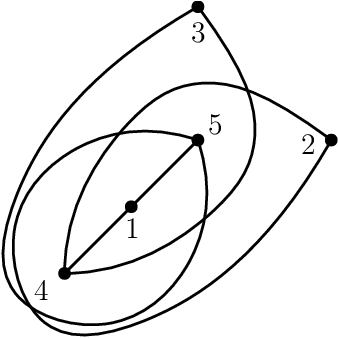}
	\caption{An almost embedding $g\colon K_{3, 2}\to \R^2$ such that $\tr(g|_{[3]*4})=\tr(g|_{[3]*5})$}
	\label{ris:K_2,3}
	\end{figure}	
	\item\label{r,e:apr-k5}  
	By \cite[Proposition~5.1]{AMS} it follows that
	%\begin{equation}
$$l(g)=-w_g(234,1)+w_g(134,2)-w_g(124,3)+w_g(123,4),$$
%\end{equation}
	So one might try to prove Theorem~\ref{t:main} by ignoring the given almost embedding of the graph $\K$ and considering only its restriction to $K_4$.
	But \cite[Figure~4.2]{AMS} shows that this approach does not work.

	\item\label{r,e:apr-k3} 
	By \cite[Proposition~5.1]{AMS} and the equality in \cite[Remark~8.4.b]{AMS} it follows that
	%\begin{equation}%\label{e:l_k13} 
	$$l(g)=\tr(g|_{[3]*4})-\tr(g|_{[3]*5}). $$
	%\end{equation} 
	So one might try to prove Theorem~\ref{t:main} by ignoring the given almost embedding of the graph $\K$ and considering only its restriction to $K_{3,2}$.
	But Figure~\ref{ris:K_2,3} shows that this approach does not work either (a more complicated example with $\tr(f|_{[3]*4})-\tr(f|_{[3]*5})\ne\pm2$ was found earlier by A.~Lazarev). 
	
	\item\label{r,e:apr-sep}
	To prove Theorem~\ref{t:main}, it suffices to prove it for a PL almost embedding in general position, see the definition before Remark~\ref{r:vaKa}.
	The natural stronger assertion that for any PL almost embedding \(g\colon \K \to \mathbb{R}^2\) in general position the points \(g(4)\) and \(g(5)\) lie in adjacent components of \(\mathbb{R}^2 \setminus g(K_3)\) is false, see Figure~\ref{ris:non-adj}.

	\begin{figure}[h]
		\centering
		\includegraphics[scale=0.4]{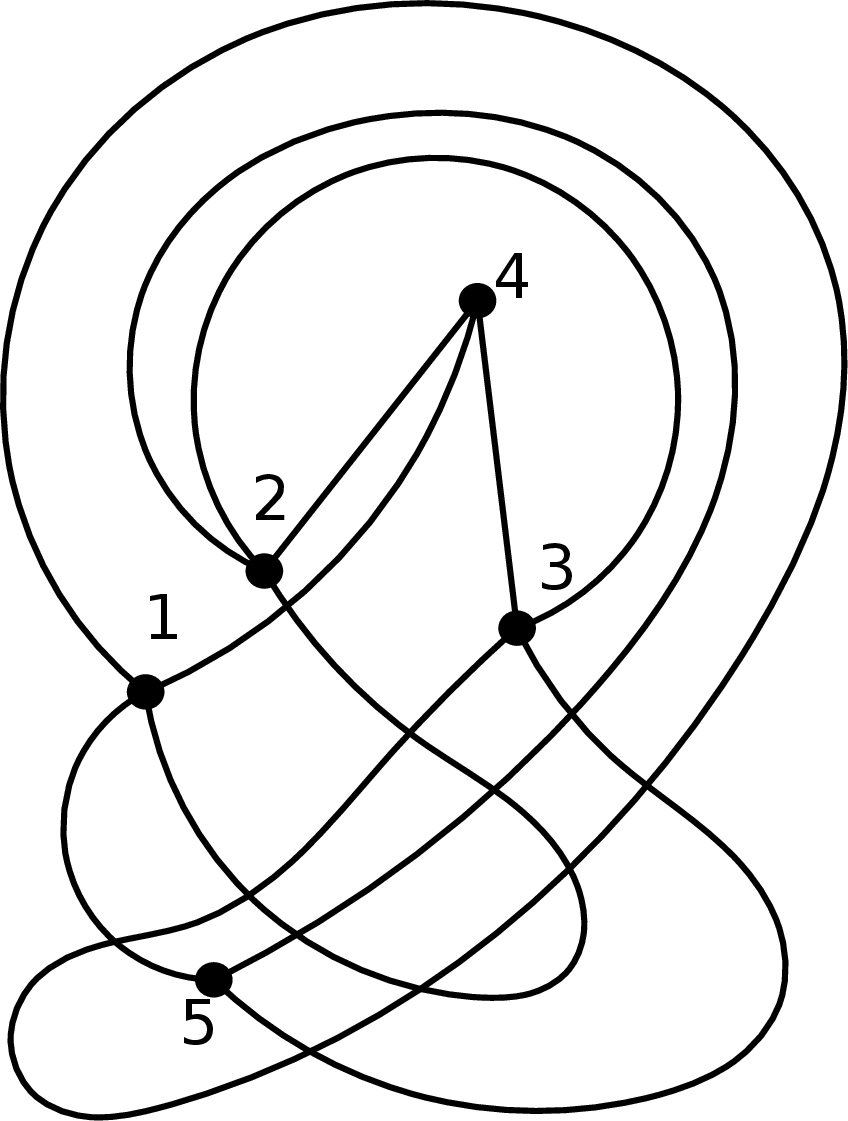}
		\caption{An almost embedding $g:\K\to\R^2$ such that \(g(4)\) and \(g(5)\) lie in non-adjacent components.}
	\label{ris:non-adj}
	\end{figure}
	
	\item\label{r,e:apr-sec_bul}	
	The assumption in the first bullet of Theorem~\ref{t:sim-form} is analogous to the condition $g(i)\notin g(jk)$ in Proposition~\ref{p:low-dim}.
	One might try to prove Theorem~\ref{t:sim-form} by ignoring the second bullet. 
	Figure~\ref{ris:adi_con} shows that the assumption in the second bullet of Theorem~\ref{t:sim-form} is also necessary.
	
	\begin{figure}[h]
		\centering
		\includegraphics[scale=1]{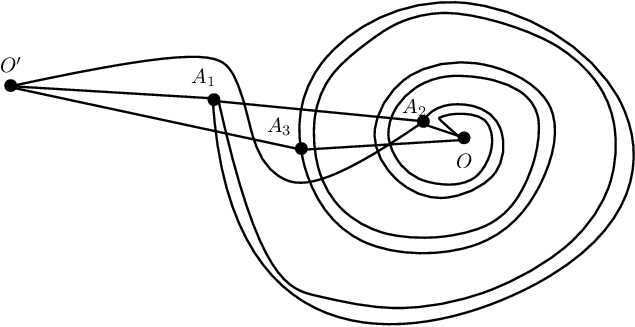}
		\caption{A counterexample to Theorem~\ref{t:sim-form} without the assumption in the second bullet: here $w(l,O)-w(l,O')=-3$.}
	\label{ris:adi_con}
	\end{figure}
	\end{enumerate}
\end{remark}

By equation in Remark~\ref{r,e:apr-k5}, \cite[Figure~4.2]{AMS}, Theorem~\ref{t:main} and by equation in Remark~\ref{r,e:apr-k5}, figure~\ref{ris:K_2,3}, Theorem~\ref{t:main} the following corollary holds.
	
\begin{corollary}\label{c:non-ext}
There are almost embeddings $K_{3,2}\to\R^2$ and $K_4\to\R^2$ none of which is extendable to an almost embedding  $\K\to\R^2$.
\end{corollary}

The analogues of Corollary~\ref{c:non-ext} for embeddings are false. 
(Every embedding $K_{3,2}\to\R^2$ and every embedding $K_4\to\R^2$ extends to an embedding $\K\to\R^2$.)

\begin{remark}[more general context: Wu invariant of (almost) embeddings]\label{r:wu}
For any graph $K$ denote 
$$\widetilde{K}:=\cup \{\sigma\times \tau  \subset K\times K: \sigma, \tau  \text{ are non-adjacent edges of } K\}.$$ 
Then the map $\widetilde g\colon \widetilde{K} \to S^1$ is
well-defined by the Gauss formula $\widetilde g(x, y):=\dfrac{g(x)-g(y)}{|g(x)-g(y)|}$. 
Define the involution $t\colon \widetilde{\K}\to \widetilde{\K}$ by $t(x, y) = (y, x)$. 
The map $\widetilde g\colon \widetilde{K} \to S^1$ is equivariant with respect to $t$ and the antipodal map of $S^1$.
The equivariant homotopy class $\alpha (g)$ of the map $\widetilde g$ is called \textit{the Haefliger-Wu invariant} of $g$ \cite[\S5]{Sk06}. 
The integers from Remark \ref{r:all_odd} are `parts' of  $\alpha(g)$. 
For definitions of `cohomological' version of this invariant, called Wu invariant, see \cite[\S 1.6]{Sk}, %\cite[\S 1.6]{Sk24},
 \cite[Theorem 4.4]{Sk06}.  

Theorem \ref{t:main} gives a restriction on values of the Haefliger-Wu invariant for almost embeddings $\K\to \R^2$. 
So Theorem \ref{t:main} is a (presumably important) step towards the interesting problem of describing the values of $\alpha (g)$ for almost embeddings $g$ of an arbitrary graph \cite[Open Problem 7]{AM25}.
\end{remark}

Some points in the plane are in \textit{general position} if no three of them lie in a line and no three segments joining them have a common interior point.

A PL map $g\colon K \rightarrow \R^2$ is in \emph{general position}, if there is a  subdivision $K'$ of $K$ such that the corresponding map $g'\colon K'\to \R^2$ is linear on any edge and the images of vertices in $K'$ in general position.

\begin{remark}[The van Kampen-Flores theorem and the lack of its integer version]\label{r:vaKa}
\begin{enumerate}[(a)]
\item\label{r,e:vaKa-def} For a general position PL map $g\colon K_5\to \R^2$ let the \emph{van Kampen number} $v(g)\in \mathbb{Z}_2$ be the sum mod $2$ of the numbers $|g\sigma \cap g\tau|$ over all non-ordered pairs $\{\sigma, \tau\}$ of non-adjacent edges of $K_5$.
The van Kampen-Flores theorem for the plane states that this number is odd, see survey \cite[Lemma 1.4.3]{Sk18}.
This is equivalent to $l(g)$ being odd, see proof in \ref{r,r:vaKa-eq}.

Theorem \ref{t:main} is an integer analog of $l(g)$ being odd. 
However, it is known that there are no integer analogs of the van Kampen-Flores theorem, see a rigorous formulation in \ref{r,e:vaKa-int_ana_ref}.

\item\label{r,r:vaKa-eq} 
%The equivalence holds because 
For any PL general position almost embedding $g\colon  \K\to \R^2$ and some general position extension $\bar{g}\colon K_5\to \R^2$ of $g$ we have $v(\bar{g})\underset2\equiv l(g)$ because  
$$v(\bar{g})\underset2\equiv \sum_{\{\sigma, \tau\}} |\bar{g}\sigma\cap \bar{g}\tau| \underset2\equiv|\bar{g}(45)\cap \bar{g}(12)|+|\bar{g}(45)\cap \bar{g}(23)|+|\bar{g}(45)\cap \bar{g}(31)| \stackrel{(1)}{\underset2\equiv} l(g),$$ where congruence $(1)$ follows from \cite[Statement~4.2.(a)]{AMS1}. 

\item\label{r,e:vaKa-int_ana_ref} Proposition. {\it Take some orientations on edges in $K_5$.
For any cell subcomplex $C$ of $\widetilde{K_5}$ and PL general position map $g\colon K_5\to \R^2$ denote  $V_C(g):=\sum_{\sigma\times \tau\subset C} g(\sigma) \cdot g(\tau)$, where $\sigma, \tau$ are edges in $K_5$, and $g(\sigma) \cdot g(\tau)$ is the sum of signs of intersection points of $g(\sigma)$ and $g(\tau)$ (for definition of $\widetilde{K}$ see Remark~\ref{r:wu}).
If $V_C$ does not depend on such $g$, then $V_C(g)=0$.  }

\emph{Proof.} 
Define the involution $t\colon \widetilde{K_5}\to \widetilde{K_5}$ by $t(x, y) = (y, x)$. 
Since $g(\sigma) \cdot g(\tau)=-g(\tau) \cdot g(\sigma)$, we may assume that $C$ and $t(C)$ have no common $2$-cells. 
Since $V_C$ does not depend on $g$, by \cite[Lemma 3.3]{Sh57} the projection of $C$ to $\widetilde{K_5}/t$ is an integer $2$-cycle for some orientation on 2-cells of $\widetilde{K_5}/t$.
So $C=\varnothing$ by known fact \ref{r,e:vaKa-pro_surf}.  

\item\label{r,e:vaKa-pro_surf} Proposition. \emph{The 2-complex $\widetilde{K_5}/t$ has only empty integer $2$-cycle (for definition of $\widetilde{K}$ and $t$ see Remark~\ref{r:wu}).}

\emph{Proof.}
It suffices to show that $\widetilde{K_5}/t$ is a compact connected non-orientable $2$-manifold, because then $\widetilde{K_5}/t$ has only empty integer $2$-cycle, see e.g. \cite[6.2, second bullet]{Mu}, cf. \cite[3.4]{Sa91}.

Since $\widetilde{K_5}/t$ is a finite cell complex, $\widetilde{K_5}/t$ is a compact complex.

Since the link of any vertex of $\widetilde{K_5}/t$ is circular, $\widetilde{K_5}/t$ is a 2-manifold.

The connectivity of $\widetilde{K_5}/t$ is proved by checking that for each two vertices of $\widetilde{K_5}/t$ there is a cell, 
containing these vertices.

By counting cells we see that $\widetilde{K_5}$ has Euler characteristic $-10$. 
For any cell $\sigma \subset \widetilde{K_5}$ the preimage $t^{-1}(\sigma)$ consists of two cells of the same dimension.
Hence $\widetilde{K_5}/t$ has Euler characteristic $\frac{-10}{2}=-5$.
Since the Euler characteristic of a closed orientable $2$-manifold is even, see e.g. \cite[6.1, twelfth bullet]{Mu}, the 2-manifold $\widetilde{K_5}/t$ is non-orientable. 

\item\label{r,e:vaKa-kam_k33} For a general position PL map $g\colon K_{3,3}\to \R^2$ let the van Kampen number $v(g)\in \mathbb{Z}_2$ be the sum mod $2$ of the numbers $|g(\sigma) \cap g(\tau)|$ over all non-ordered pairs $\{\sigma, \tau\}$ of non-adjacent edges of $K_{3,3}$.
It is known that this number is odd, see survey \cite[Remark 1.4.4.a]{Sk18}.
Analogously to \ref{r,e:vaKa-def}, \ref{r,r:vaKa-eq} the number $l(g)$ and $v(g)$ are congruent modulo $2$, and there are no integer analogs of the van Kampen number.

\end{enumerate}
\end{remark}

\begin{lemma}\label{p:GreatLemma2}
	For  any almost  embedding $g\colon \K \rightarrow \R^2$ there is a PL almost embedding $f\colon \K \rightarrow \R^2$ such that $f|_{K_{3,2}}$ is an embedding and $l(f) = l(g)$.
\end{lemma}

We illustrate our main idea by deducing Theorem~\ref{t:main} from Lemma~\ref{p:GreatLemma2}.
In that proof we reduce Theorem~\ref{t:main} to a version of Proposition~\ref{p:low-dim}.  
%Proposition~\ref{p:idea-el} is equivalent to Proposition~\ref{p:low-dim}.

\begin{proof}[Proof of Theorem \ref{t:main} assuming Lemma \ref{p:GreatLemma2}]
	By Lemma \ref{p:GreatLemma2} it suffices to prove Theorem~\ref{t:main} under the additional assumption that $g|_{K_{3,2}}$ is a PL embedding. 
	Recall the following known result. 
	
	\emph{
	For any embeddings $f, g:K_{3, 2}\to \R^2$ there is a homeomorphism $s:\R^2\to\R^2$ such that $f(K_{3,2})=s(g(K_{3,2}))$.
	}

	Therefore, it suffices to prove Theorem~\ref{t:main} under the additional assumption that $g(K_{3,2})$ coincides with the image of the standard embedding.
	
	Indeed, the statement above implies that there exists a homeomorphism
	$s\colon\R^2\to\R^2$ such that $s(g(K_{3,2}))$ coincides with the image of the standard embedding of $K_{3,2}$ in $\R^2$.
	The homeomorphism $s$ preserves winding numbers up to sign.
	Hence $$l(g) = w_g(123, 4)-w_g(123, 5) = \pm(w_{s\circ g}(123, 4)-w_{s\circ g}(123, 5)) = \pm l(s\circ g).$$
	
	Now assume that $g(K_{3,2})$ coincides with the image of the standard embedding.
	Under this assumption, Theorem~\ref{t:main} is equivalent to the following version of Proposition~\ref{p:low-dim}.
		
	\textit{ Let $A_1, A_2, A_3$ be the vertices of a regular triangle in the plane and $O$ its center.
			For $m\in [3]$ %аналогично
			let $l_m$ be a polygonal line  joining two of the vertices distinct from $A_m$ and disjoint from the ray $OA_m$.
			Orient the lines so that $l_1l_2l_3$ is a closed curve.
			Then $w(l_1l_2l_3, O) = \pm 1$.}
	 
\end{proof}

\begin{remark}\label{r:just}
\begin{enumerate}[(a)]
\item\label{r,e:just-geo} 

	Since the proof of Lemma~\ref{p:GreatLemma2} (and hence of Theorem~\ref{t:main}) is technical, it is natural to look for a short proof of Theorem~\ref{t:main}  using algebraic topology. 
	
	One possible approach is to apply Borsuk–Ulam-type theorems, such as the van Kampen-Flores theorem, see Remark~\ref{r:vaKa}.\ref{r,e:vaKa-def}. 
	However, this yields only a restriction on the parity, rather than the stronger restriction of the values be $\pm1$.

	Another approach is to express \(l(g)\) using simpler invariants such as homological relations in \(\widetilde{K}\) (see Remark~\ref{r:apr}.\ref{r,e:apr-k5} and \ref{r:apr}.\ref{r,e:apr-k3}). 
	However, this expression does not make Theorem~\ref{t:main} trivial.

	%One may try to apply the Borsuk--Ulam likes theorems. 
	%Tthis aproach lead only to oddness, but not in terms of $\pm 1$, see Remark~\ref{r:vaKa}.\ref{r,e:vaKa-def}.
	
	%One may try to express $l(g)$ using homological relations in $\widetilde{K}$ (see Remark~\ref{r:wu}). 
	%This approach expresses $l(g)$ in terms of simpler invariants, but this alone does not make Theorem~\ref{t:main} trivial.

	For these reasons the restriction   of the form $\pm 1$ in Theorems~\ref{t:sim-form} and~\ref{t:main} motivated us to use geometric rather than algebraic aproach. 
	(Proposition~\ref{p:low-dim} also contains restriction of the form $\pm 1$, but this statement is so simple that the distinction between an algebraic and a geometric proof is blurred.)
 
	%This suggests that in Theorem~\ref{t:main} one should use geometric techniques.

\item\label{r,e:just-tech}	
	Some theorems in PL topology of the plane have technical proofs, while attempts for simpler proofs led to mistakes. 
	An example is the completeness of the van Kampen planarity  obstruction.
	%For a description and recovery of a gap in \cite[\S2]{Sa91} see  \cite[footnote 6???10]{Sc13} and \cite[proof of Theorem 3.1]{MPS}.
	For a mistake in the proof of \cite[(2.2)]{Sa91} see \cite[footnote 6]{Sc13}; for a correct proof see e.g. \cite[proof of Theorem 3.1]{MPS}.
	
	This justifies the need of a careful proof of Theorem \ref{t:main} (or of Lemma \ref{p:GreatLemma2}), and explains why that proof is technical. 
\end{enumerate}
	
	%Раз в теореме 1.1 утверждение скорее геометрическое, а не алгебраическое, то доказываться оно должно геометрическими методами, а не алгебраическими.
	
	%The existence of the van Kampen obstruction can be proved by algebraic methods, but the obstruction in that case concerns oddness. 
	
\end{remark}

\section{Lemmas for proof of Lemma~\ref{p:GreatLemma2}}\label{s:lemmas}

A PL almost  embedding $f\colon \K \rightarrow \R^2$ in general position is called an {\bf improvement} of PL almost  embedding $g\colon \K \rightarrow \R^2$, if $l(f) = l(g)$ and the number of the self-intersection points of $f$ is fewer than the number of the self-intersection points of $g$.
%A PL almost  embedding $g\colon \K \rightarrow \R^2$ is {\bf interesting} if there is an improvement of $g$.

The following Lemma~\ref{l:LessLemma3} is the main tool in the proof of Lemma~\ref{p:GreatLemma2}, enabling the elimination of self‑intersections.

In this section, except in the proof of Lemma~\ref{p:GreatLemma2}, we assume that $g\colon \K\to\R^2$ is a PL almost embedding in general position.
In this section any arc is a part of some edge of $\K$.

We write $gA$ instead of $g(A)$. 

\begin{lemma}[a version of Whitney trick]\label{l:LessLemma3}
    Assume that for arcs $I, \bar{I}\subset K_{3,2}$ we have
        
    	\begingroup\makeatletter
    	\edef\@currentlabel{\getrefnumber{l:LessLemma3}.0}%
    	\phantomsection\label{l:LessLemma3.0}%
    	\endgroup
    	(\getrefnumber{l:LessLemma3}.0)\;	 
    	$g|_{I}$ and $g|_{\bar{I}}$ are embeddings; 
        
        \begingroup\makeatletter
    	\edef\@currentlabel{\getrefnumber{l:LessLemma3}.1}%
    	\phantomsection\label{l:LessLemma3.1}%
    	\endgroup
    	(\getrefnumber{l:LessLemma3}.1)\;
        $gI\cap g\bar{I}=g\partial I=g\partial \bar{I}$;
                
        \begingroup\makeatletter
    	\edef\@currentlabel{\getrefnumber{l:LessLemma3}.2}%
    	\phantomsection\label{l:LessLemma3.2}%
    	\endgroup
    	(\getrefnumber{l:LessLemma3}.2)\;
         $g[5]$ is contained in the closure of one of the two connected component of $\R^2\backslash (gI\cup g\bar{I});$ 
        
    Then there is an improvement of $g$.
 
\end{lemma}

\begin{remark}\label{r:addendum}
	Assume that some intervals $I, \bar{I}$  satisfy the properties (\ref{l:LessLemma3.0}) and (\ref{l:LessLemma3.1}).
	Then $\R^2\backslash (g(I)\cup g(\bar{I}))$ has exactly two connected components, say $U$ and $W$.
	Additionally, we have  $g(I)\cup g(\bar{I}) = \partial U = \partial W$.
	Hence for any $pt\in I$ we have $g(pt)\in \partial U = \partial W$.
\end{remark}

For $g\colon X\to Y$ and for $A\subset Y$ denote $g^{-1}A:=\{x\in X\ :\ g(x)\in A\}$.

For any points $p, q$ in edge $ij$ of the graph $\K$ denote by $[p, q]\subset ij$ the part of the edge $ij$ between $p$ and $q$.   

\begin{proof}[Proof of Lemma~\ref{l:LessLemma3}]
    Denote by $U$ the complement to the closure of the connected component from the property (\ref{l:LessLemma3.2}).
    For any arc $J$ denote by $\alpha_J$ the edge of $\K$ containing $J$.
    For any arc $J \subset g^{-1}(\Cl U)$ denote by $J^{\varepsilon}$ the intersection of some small open neighborhood of $J$ and $\alpha_J$.
    
    {\it Case 1: there is an arc $J$ in $g^{-1}(\Cl U)$ such that $g|_J$ is not embedding.}  
    Modify  $g$ on $J$ to match $f$ as shown in Figure~\ref{ris:Reid} such that $fJ \subset \Cl U$, and such that for any edge $\beta$ distinct from $\alpha_J$ we have $g\beta\cap g\alpha_J=f\beta\cap f\alpha_J$.
    Hence $f$ is almost embedding.
    The number of self-intersection points of $f$ is fewer than the number of self-intersection points of $g$.
    
    \begin{figure}[h]
    \centering
        \includegraphics[scale=0.6]{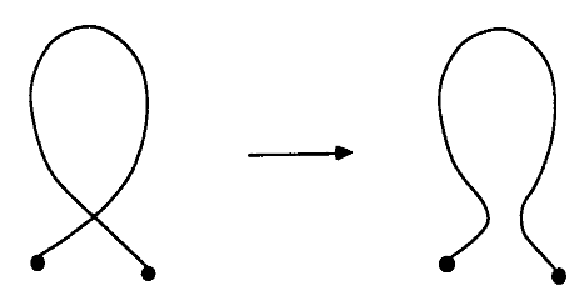}
         \caption{Transformation of a map}
         \label{ris:Reid}
\end{figure}
    
    {\it Proof that $l(f) = l(g)$ in Case 1.}  
    If $J \not\subset 123$, then the restrictions of $f$ and $g$ to $123$ and $\{4, 5\}$ coincide. 
    So $l(f) = l(g)$. 
    Assume that $J \subset 123$.
    Let $\Gamma \colon  S^1 \to \mathbb{R}^2$ be a map such that 
    
    $\bullet$ $\Gamma$ is a map onto $gJ\cup fJ$;
    
    $\bullet$ $\Gamma$-image of one semicircle is $gJ$ and the image of the other semicircle is $fJ$;
    
    $\bullet$ $\Gamma \{-1, 1\}=g\{\partial J\}$;
    
    $\bullet$ for any $s\in \{4, 5\}$ we have $w_f(123, s)=w_g(123, s)+ w(\Gamma, gs)$. 
    
    Since $\Gamma S^1\subset \Cl U$, and $g4, g5 \notin U$, and there is a curve $p\colon  [0, 1] \to \R^2 \setminus U$ joining points $g4$ and $g5$ such that $p(\Int [0, 1] ) \cap \Cl U = \varnothing$, we have $w(\Gamma, g4)=w(\Gamma, g5)$. 
    Hence
    $$l(f) = w_f(123, 4)- w_f(123, 5)=  (w_g(123, 4)+ w(\Gamma, g4))- (w_g(123, 5)+ w(\Gamma, g5))=$$
    $$ = w_g(123, 4)- w_g(123, 5) + (w(\Gamma, g4)-w(\Gamma, g5))=w_g(123, 4)- w_g(123, 5)=l(g).$$
    
    Then $f$ is improvement of $g$.
    
    {\it Case 2: for any arc $J$ in $g^{-1}(\Cl U)$ the restriction $g|_J$ is an embedding.}
    An arc $J\subset g^{-1}(\Cl U)$ is called \textbf{extendable} if there is an arc $J'\subset g^{-1}(\Cl U)$ such that $J\neq J'$ and $g(\partial J)=g(\partial J')$.
    The arc $\bar J\subset g^{-1}(\Cl U)$ denotes some of this counterpart of $J$ satisfying these conditions.
    Since $g$ is in general position and $g$ restricted to any arc in $ g^{-1}(\Cl U)$ is embedding, we have that for any extendable arc $J$ there is a unique $\bar{J}$.
    For an extendable arc $J$ denote by $U_J$ the open component in $\R^2\backslash g(J\cup \bar{J})$ such that $U_J\subset U$.
    An extendable arc $J$ is called \textbf{minimal} if for every extendable arc $J'\notin \{J, \bar J\}$ we have $U_{J'}\not\subset U_J$.
    
    Let $J$ and $J'\notin \{J, \bar{J}\}$ be extendable arcs. 
    Since $\partial U_J=g(J\cup \bar{J})\neq g(J'\cup \bar{J'})=\partial U_{J'}$, we have $U_J\neq U_{J'}$. 
    Hence if $J$ and $J'\notin \{J, \bar{J}\}$ are extendable arcs such that $U_{J'} \subset U_J$, then $U_{J'} \varsubsetneq U_J$.

    Let us show that there is a minimal arc $J$.
    It follows from properties~(\ref{l:LessLemma3.0}) and~(\ref{l:LessLemma3.1}) that $I$ is an extendable arc.
    Denote $J_0:=I$.
    If $J_0$ is minimal, then we are done.
    Otherwise there is an extendable arc $J_1\notin \{J_0, \bar{J_0}\}$ such that $U_{J_1}\varsubsetneq U_{J_0}$.
    If $J_1$ is minimal, then we are done.
    Otherwise there is an extendable arc $J_2\notin \{J_1, \bar{J_1}\}$ such that $U_{J_2}\varsubsetneq U_{J_1}$.
    Proceeding in this way, we produce a sequence $$U_{J_0}\varsupsetneq U_{J_1}\varsupsetneq U_{J_2}\varsupsetneq \ldots$$ 
    Since the set of extendable arcs is finite, we have that the sequence must terminate.
    Hence there is a minimal arc $J$. 
     
    %Denote $\gamma_J:=J\cup \bar{J}$.
    Take  $J^{\varepsilon}, \bar{J}^{\varepsilon}$ sufficiently small such that for every edge $\beta$ distinct from $\alpha_{J^{\varepsilon}}$ and $\alpha_{\bar{J}^{\varepsilon}}$ we have 
    \begin{equation}
    g\beta\cap g(J^{\varepsilon}\cup \bar{J}^{\varepsilon})=g\beta\cap g(J\cup \bar{J}).
    \tag{*}
    \end{equation}
    Take sufficiently small open neighbourhood $U_J^{\varepsilon}$ of $U_J$ such that $J^{\varepsilon}, \bar{J}^{\varepsilon}\subset g^{-1}\Cl U_J^{\varepsilon}$, and $g4, g5\notin U_J^{\varepsilon}$, and there is a curve $p\colon  [0, 1] \to \R^2 \setminus U_J^{\varepsilon}$ joining points $g4$ and $g5$ such that $p(\Int [0, 1] ) \cap \Cl U_J^{\varepsilon} = \varnothing$.
    Modify  $g$ on ${J^{\varepsilon}\cup \bar{J}^{\varepsilon}}$ to match $f$ as shown in Figure~\ref{ris:Lemm3.1} such that $f(J^{\varepsilon}\cup \bar{J}^{\varepsilon})\subset \Cl U_J^{\varepsilon}$, and such that for any edge $\beta$ distinct from $\alpha_J$ and $\alpha_{\bar{J}}$ we have 
    \begin{equation}
    g\beta\cap f(J^{\varepsilon}\cup \bar{J}^{\varepsilon})=g\beta\cap g(J^{\varepsilon}\cup \bar{J}^{\varepsilon}).
    \tag{**}
    \end{equation}
    \begin{figure}[h]
    \centering
\includegraphics[scale=0.8]{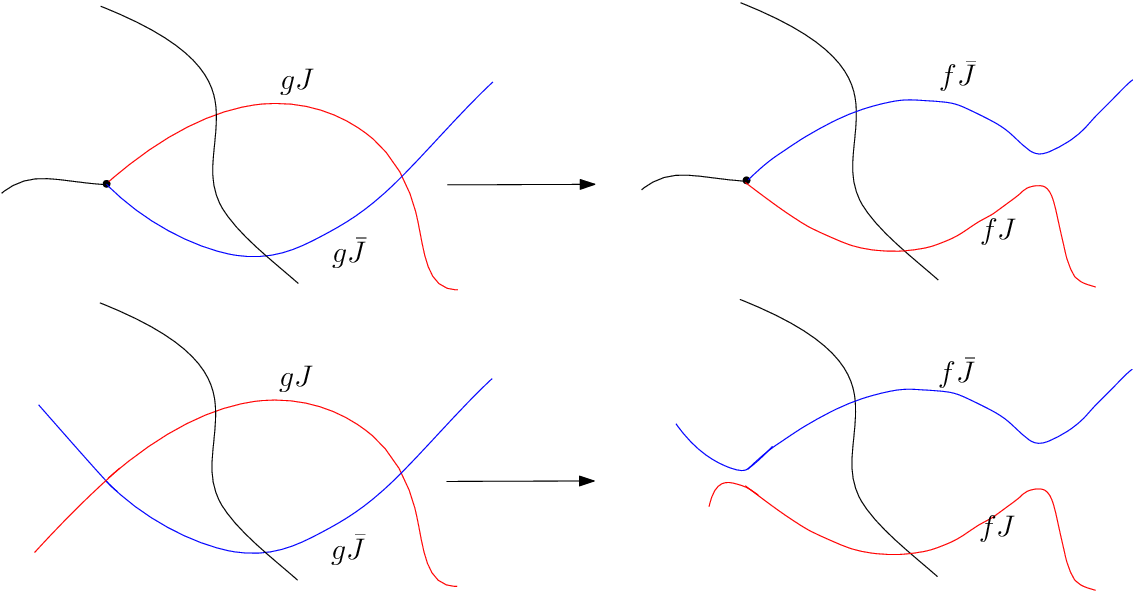}
	 \caption{The upper figure is used if $\partial J \cap [5]\neq\varnothing$ or $\partial \bar{J} \cap [5]\neq\varnothing$, and the lower figure is used otherwise}
    \label{ris:Lemm3.1}
\end{figure}

    \textit{Proof that $l(f) = l(g)$ in Case 2.}
    If $J$ and $\bar{J}$ are not in $123$, then the restrictions of $f$ and $g$ coincide on $123$ and $\{4, 5\}$.
    So $l(f) = l(g)$. 
    
    The proof that $l(f) = l(g)$ under the additional assumption that $J$ is a subset of $123$ is obtained from the proof of $l(f) = l(g)$ under the additional assumption that $J\subset 123$ in the case 1 by replacing $J$ with $\Cl J^{\varepsilon}$ and $U$ with $U_J^{\varepsilon}$.
    
    Analogously $l(f) = l(g)$ under the additional assumption that $\bar{J}$ is subset of $123$.
    
    Assume that $J, \bar{J}\subset 123$ (cf. the penultimate paragraph in case 1).
    Let $\Gamma, \bar{\Gamma} \colon  S^1 \to \mathbb{R}^2$ be a map such that 
    
    $\bullet$ $\Gamma$ is a map onto $g\Cl J^{\varepsilon}\cup f\Cl J^{\varepsilon}$ and $\bar{\Gamma}$ is a map onto $g\Cl \bar{J}^{\varepsilon}\cup f\Cl \bar{J}^{\varepsilon}$;
    
    $\bullet$ $\Gamma$-image of one semicircle is $g\Cl J^{\varepsilon}$ and the $\Gamma$-image of the other semicircle is $f\Cl J^{\varepsilon}$, and $\bar{\Gamma}$-image of one semicircle is $g\Cl \bar{J}^{\varepsilon}$ and the $\bar{\Gamma}$-image of the other semicircle is $f\Cl \bar{J}^{\varepsilon}$;
    
    $\bullet$ $\Gamma \{-1, 1\}=g\{\partial J^{\varepsilon}\}$ and $\bar{\Gamma} \{-1, 1\}=g\{\partial \bar{J}^{\varepsilon}\}$;
    
    $\bullet$ for any $s\in \{4, 5\}$ we have $w_f(123, s)=w_g(123, s)+ w(\Gamma, gs)+ w(\bar{\Gamma}, gs)$.
    
    By the choose of $U_J^{\varepsilon}$ and by $\Gamma S^1, \bar{\Gamma}S^1\subset \Cl U_J^{\varepsilon}$, we have $w(\Gamma, g4)=w(\Gamma, g5)$ and $w(\bar{\Gamma}, g4)=w(\bar{\Gamma}, g5)$. 
    Hence
    $$l(f) = w_f(123, 4)- w_f(123, 5)=  (w_g(123, 4)+ w(\Gamma, g4)+w(\bar{\Gamma}, g4))- (w_g(123, 5)+ w(\Gamma, g5)+w(\bar{\Gamma}, g5))=$$
    $$ = w_g(123, 4)- w_g(123, 5) + (w(\Gamma, g4)-w(\Gamma, g5) + w(\bar{\Gamma}, g4) - w(\bar{\Gamma}, g5))=w_g(123, 4)- w_g(123, 5)=l(g).$$
    
    \textit{Proof that $f$ is almost embedding in Case 2.}
    The restrictions of $f$ and $g$ to the complement of $J^{\varepsilon}\cup \bar{J}^{\varepsilon}$ in $\K$ coincide. 
    Hence for any non-adjacent edges $\beta, \beta'$ distinct from $\alpha_{J}$ and $\alpha_{\bar{J}}$ we have $$f\beta \cap f\beta' = g\beta \cap g\beta' = \varnothing.$$ 
    
    Consider any edge $\beta$ of $\K$ non-adjacent to  $\alpha_{J}$. 
    Since $g$ is an almost embedding, we have $g\beta \cap g\alpha_{J}=\varnothing$. 
    Since $g$ restricted to any arc in $g^{-1}(\Cl U)$ is embedding and $g\beta \cap g\alpha_{J}=\varnothing$, we have that $g\beta\cap \Cl U_J$ consists of images of extendable arcs $J'$ such that $g(\partial J')\subset g\bar{J}$. 
    Hence $U_{J'}\subset U_{\bar{J}}=U_J$.
    Since $J$ is a minimal arc, we have $g\beta \cap g(J\cup \bar{J})\subset g\beta\cap\Cl U_J=\varnothing$.
    We have $$f\beta\cap f\alpha_{J}=(f\beta\cap f(\alpha_{J}\backslash J^{\varepsilon}))\cup(f\beta\cap fJ^{\varepsilon})\stackrel{(1)}{\subset} (g\beta\cap g(\alpha_{J}\backslash J^{\varepsilon}))\cup(g\beta\cap g(J\cup \bar{J}))\stackrel{(2)}{=}\varnothing,\quad \text{where}$$
    
    $\bullet$ inclusion $(1)$ holds because the restrictions of $f$ and $g$ to the complement of $J^{\varepsilon}\cup \bar{J}^{\varepsilon}$ in $\K$ coincide and because $f\beta\cap fJ^{\varepsilon}=g\beta\cap fJ^{\varepsilon}\subset g\beta\cap f(J^{\varepsilon}\cup \bar{J}^{\varepsilon})\overset{(^{**})}{=}g\beta\cap g(J^{\varepsilon}\cup \bar{J}^{\varepsilon})\overset{(^*)}{=}g\beta \cap g(J\cup \bar{J})$, and
    
    $\bullet$ equation $(2)$ holds because $g$ is an almost embedding and because $g\beta\cap g(J\cup \bar{J}) =\varnothing$.
    
    Analogously $f\beta\cap f\alpha_{\bar{J}}=\varnothing$   for any edge $\beta$ non-adjacent to $\alpha_{\bar{J}}$. 
    
    Hence $f$ is an almost embedding. 
    
    Then $f$ is improvement of $g$.
\end{proof}

\begin{figure}[h]
	\centering
	\includegraphics[scale=0.7]{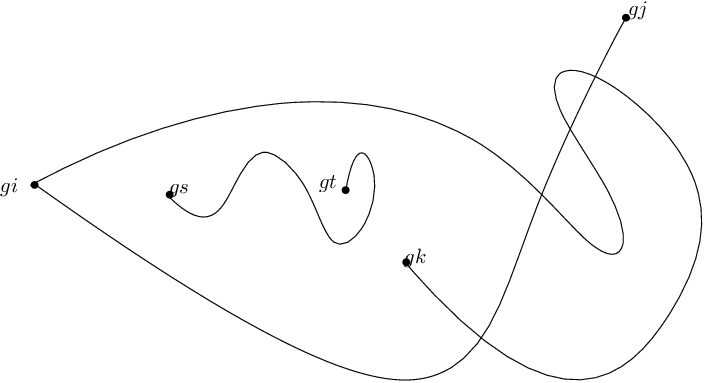}
	\caption{}
	\label{ris:LastConst}
\end{figure}
    
    We say that a triple $(i, j, k)$ of points from $[5]$ is \textbf{interesting} if 
    
    $\bullet$ $ij, ik$ are edges of $K_{3,2}$,
    
    $\bullet$ $g|_{ij}$ and $g|_{ik}$ are embeddings, and
    
    $\bullet$ $g(ij)\cap g(ik)\varsupsetneq\{gi\}$.
    
    Lemma~\ref{l:improv} shows that interesting triples provide intervals that almost satisfy the conditions of Lemma~\ref{l:LessLemma3}.
    
    For an interesting triple $(i, j, k)$ denote by $p=p_{i, j, k, g}\in ij$ the first point after $i$ on the path $g|_{ij}$ (starting from $i$) contained in $g(ik)$. 
    Denote by $q=q_{i, j, k, g}\in ik$ the point such that $gp=gq$.
    Points $p, q$ are called \textbf{improving for the interesting triple} $(i, j, k)$.
    An open component $U$ in $\R^2\backslash (g[i,p]\cup g[i, q])$ is called \textbf{improving for the interesting triple} $(i, j, k)$ if  $g([5]\backslash \{i, j, k\})\subset~U$. 
	
    A triple $(i, j, k)$ is interesting iff $(i, k, j)$ is interesting. 
    The orders of the points $j$ and $k$ is important for the definition of improving points and improving component.
	
	\begin{lemma}\label{l:improv}
		Assume that $(i, j, k)$ is an interesting triple.
		Denote by $p\in ij$ and $q\in ik$ the improving points for the interesting triple $(i, j, k)$.
		We have that:
		
		(1) arcs $I:=[i, p]$ and $\bar{I}:=[i, q]$ lie in $K_{3, 2}$ and satisfy the properties (\ref{l:LessLemma3.0}), (\ref{l:LessLemma3.1}); 
		
		(2) there is a unique open component $U$ such that $U$ is improving for $(i, j, k)$, see Figure \ref{ris:LastConst}.
	\end{lemma}

	\begin{proof}	
	\textit{Proof of (1).} 
	By the definition of the interesting triple $I, \bar{I} \subset K_{3,2}$. 
	Since the restriction of $g$ to any edge of $K_{3,2}$ is an embedding and $ij, ik$ are edges of $K_{3,2}$, we have that $g|_{I}$ and $g|_{\bar{I}}$ are embeddings.
	Hence arcs $I$ and $\bar{I}$ satisfy the property (\ref{l:LessLemma3.0}). 
	Since $p$ is the first point after $i$ on the path $g|_{ij}$ (starting from $i$) contained in $g(ik)$, we have $g[i, p]\cap g[i, q]=\{gi, gp\}$. 
	Hence arcs $[i, p]$ and $[i, q]$ satisfy the property (\ref{l:LessLemma3.1}). 
	
	\textit{Proof of (2).} 
	Denote by $e$ the edge in $\K$ such that $\partial e \subset [5]\backslash \{i, j, k\}$.
	We have $$ge\cap (g[i,p]\cup g[i, q])\subset ge\cap (g(ij)\cup g(ik))=(ge\cap g(ij))\cup(ge\cap g(ik)) \stackrel{(*)}{=} \varnothing,$$ 
	where equation $(*)$ holds because $g$ is a PL almost embedding. 
	It follows that there is a unique component $U$ in $\R^2\backslash (g[i,p]\cup g[i, q])$ that contains $gs$ and $gt$. 	
	%In this paragraph we show that for $U$ from the previous paragraph we have $gi\in\partial U$.
	%Since $g|_{ij}$ and $g|_{ik}$ are embeddings, and $g[i,p]\cap g[i, q]=\{gi, gp\}$, we have $\R^2\backslash (g[i, p]\cup g[i, q])$ has two connected components.
	%Hence $g[i, p]\cup~g[i, q]=\partial U$.
	%Then $gi\in g[i, p]\cup~g[i, q]=\partial U$. 
\end{proof}   

\begin{lemma}\label{l:improv_exist}
	Assume that the restriction of $g$ to any edge of $K_{3,2}$ is an embedding.
	Assume that an open component $U$ is improving for the interesting triple $(i, j, k)$.
	If either $gj, gk\in U$, or $gk\notin U$, then there are arcs satisfying the conditions of Lemma~\ref{l:LessLemma3}.
\end{lemma}    

\begin{proof}
	Denote by $p\in ij$ and $q\in ik$ the improving points for the interesting triple $(i, j, k)$.
	
	{\it Suppose that $gj, gk \in U$.}
	By Lemma~\ref{l:improv} for the interesting triple $(i, j, k)$, we have that arcs $I:=[i, p]$ and $\bar{I}:=[i, q]$ lie in $K_{3, 2}$ and satisfy the properties (\ref{l:LessLemma3.0}), (\ref{l:LessLemma3.1}).
	Since $gj, gk \in U$, and by Remark~\ref{r:addendum} for $I, \bar{I}$ and $pt=i $ we have $gi\in \partial U$, and $U$ is improving for the interesting triple $(i, j, k)$, we have $g[5] \subset \Cl U$. 
	Then arcs $I, \bar{I}$ satisfy the property (\ref{l:LessLemma3.2}). 
	
	Hence $I, \bar{I}$ satisfy the conditions of Lemma~\ref{l:LessLemma3}.
	
	{\it Suppose that $gk\notin U$.}
	Denote by $s$ and $t$ points from $[5]\backslash\{i, j, k\}$ such that edge $ks$ is in $K_{3,2}$.
	Since $gs\in U$ and $gk\notin U$, we have $\varnothing \neq g(ks)\cap\partial U= g(ks) \cap (g[i, p] \cup g[i, q]) = g(ks)\cap g[i, q]$. 
	Since $k\notin [i, q]$ we have $gk\notin g(ks)\cap g[i, q]$.
	Hence $$g(ks) \cap g(ik) \supset (g(ks) \cap g[i, q]) \cup \{ gk \} \varsupsetneq \{ gk \}.$$
	Since $ks, ki$  are edges of $K_{3,2}$, $g|_{ks}$ and $g|_{ki}$ are embeddings, and $g(ks)\cap g(ik)\varsupsetneq \{gk\}$, we have $(k, s, i)$ is interesting triple.
	Denote by $p'\in ks$ and $q'\in ki$ the improving points for the interesting triple $(k, s, i)$.
	By Lemma~\ref{l:improv}, applied to the interesting triple $(k,s,i)$, there is a unique improving open component $U'$.
	In the next two paragraphs we show that $U\cap \partial U'=\varnothing$.
	
	Since  $U\cap \partial U'=U\cap (g[k, p']\cup g[k, q'])$, it suffices to show that $U\cap g[k, p']=\varnothing$ and $U\cap g[k, q']=\varnothing$.
	Let us show that $U\cap g[k, p']=\varnothing$. 
	Since $g|_{k s}$ has no self-intersection and $gk\notin U$, it suffices to show that $g[k, p']\cap \partial U\subset \{gp'\}$.
	We have $$g[k, p']\cap \partial U=g[k, p']\cap (g[i, p]\cup g[i, q])\stackrel{(1)}{=}$$$$\stackrel{(1)}{=}g[k, p']\cap g[i, q]\stackrel{(2)}{\subset} \{gp'\},\quad\text{where,} $$
	
	$\bullet$ equation (1) holds from the definition of an almost embedding, and
	
	$\bullet$ inclusion (2) holds because $g[k, p']\cap g(ki)=\{gk, gp'\}$ and because from the definition of $q$ we have $k\notin [i, q]$.
	
	Let us show that $U\cap g[k, q']=\varnothing$. 
	By the definition of $q$, $g(ki)\cap \partial U=g(ki)\cap (g[i, p]\cup g[i, q])=g[i, q]$.
	Since $g|_{ki}$ has no self-intersection, $gk\notin U$ and $g(ki)\cap \partial U = g[i, q]$, we have $g(ki)\cap U=\varnothing$.
	Hence $g[k, q']\cap U\subset g(ki)\cap U=\varnothing$.
	
	Since $U\cap \partial U'=\varnothing$, we have either $U\subset U'$, or $U \cap U' = \varnothing$.
	It follows from the definitions of $U$ and $U'$ that $gt\in U$ and $gt\in U'$.
	Hence $U\subset U'$.
	It follows from Remark~\ref{r:addendum}, applied with $I=[i,p]$, $\bar{I}=[i,q]$, and $pt=i$, that $gi\in\partial U$.
	It follows from Remark~\ref{r:addendum}, applied with $I=[k,p']$, $\bar{I}=[k,q']$, and $pt=k$, that $gk\in\partial U'$.
	Hence $g\{i, s, t\}\subset\Cl U$ and $g\{k, j, t\}\subset\Cl U'$.
	Since $g\{i, s, t\}\subset\Cl U$ and $g\{k, j, t\}\subset\Cl U'$, and $U\subset U'$, we have $g[5]\subset\Cl U\cup\Cl U'=\Cl U'$. 
	Then arcs $I:=[k, p']$ and $\bar{I}:=[k, q']$ satisfy the property (\ref{l:LessLemma3.2}).
	
	By Lemma~\ref{l:improv} for the interesting triple $(k, s, i)$, we have that arcs $I$ and $\bar{I}$ lie in $K_{3, 2}$ and satisfy the properties (\ref{l:LessLemma3.0}), (\ref{l:LessLemma3.1}).
	Hence $I, \bar{I}$ satisfy the conditions of Lemma~\ref{l:LessLemma3}.
\end{proof}

\begin{lemma}\label{l:improv_exist_last_case}
	Assume that the restriction of $g$ to any edge of $K_{3,2}$ is an embedding, and $(i, j, k)$ and $(k, i, l)$ are interesting triples, and $U, U'$ are open components improving for the interesting triples $(i, j, k)$ and $(k, i, l)$ respectively.
	If $gj\notin U$, $gk\in U$, $gi\notin U'$, and $gl\in U'$, then there are arcs satisfying the conditions of Lemma~\ref{l:LessLemma3}.
\end{lemma}   

\begin{proof}
	Denote by $p\in ij$ and $q\in ik$ the improving points for the interesting triple $(i, j, k)$.
	Denote by $p'\in ki$ and $q'\in kl$ the improving points for the interesting triple $(k, i, l)$.

	In this paragraph we show that $\Cl U'\not\subset \Cl U$.
	Since $p\neq j$ and $g$ in general position, we have $gj\notin g[i, p]\cup g[i, q]=\partial U$.
	Analogously $gi\notin \partial U'$.
	Since $gj\notin \partial U$ and $gi\notin \partial U'$, and $gj\notin U$ and $gi\notin U'$, we have $gj\notin \Cl U$ and $gi\notin \Cl U'$.
	It follows from the definition of $U'$ that $gj\in U'\subset \Cl U'$.
	Then we have $\Cl U'\not\subset \Cl U$.
	
	It follows from Remark~\ref{r:addendum} applied with $I=[i,p]$, $\bar{I}=[i,q]$, and $pt=i$, that $gi\in\partial U$.
	Since $gi\notin\Cl U'$, we have $gi\in \R^2\backslash\Cl U'$.
	Since $gi\in \partial U$, we have $gi\notin U$. 
	Since $gi\notin U$ and $gi\in \R^2\backslash\Cl U'$, we have $\R^2\backslash\Cl U'\not\subset U$.
	Hence we have $\partial U'\not\subset \Cl U$.
	
	In the next two paragraphs we show that there is a point $a\in [k, q']$ such that $ga\notin \Cl U$.
	We have
	\begin{align}
	g(ik)\cap \Cl U=(g[i, q]\cup g[q, k])\cap \Cl U=(g[i, q]\cap\Cl U)\cup (g[q, k]\cap\Cl U)=g(ik). \tag{*}
	\end{align}
	Let us show the last equation.
	Since $g[i, q]\subset \partial U$, we have $g[i, q]\cap\Cl U=g[i, q]$.
	Since $$g[q, k]\cap\partial U=g[q, k]\cap (g[i, q]\cup g[i, p])=(g[q, k]\cap g[i, q])\cup (g[q, k]\cap g[i, p])=gq,$$
	where the last equation holds because: 
	
	$\bullet$ since $g|_{K_{2, 3}}$ is embedding, we have $g[q, k]\cap g[i, q]=gq$, and;
	
	$\bullet$ It follows from the definition of $p$ that $g[q,k]\cap g[i,p]\subset g(ik)\cap g[i,p]=\{gi,gp\}.$
	Since $gi\notin g[q, k]$ we have $g[q, k]\cap g[i, p]\subset \{gp\}=\{gq\}$.
	
	Since $[k, p']\subset ik$ and the equation $(*)$, we have $g[k, p']\cap \Cl U=g[k, p']$.
	Since $\partial U'\not\subset \Cl U$ and $g[k, p']\cap \Cl U=g[k, p']$, we have that there is a point $a\in [k, q']$ such that $ga\notin \Cl U$.
	
	Since $gk, gl\in \Cl U$, $ga\notin \Cl U$ and $\{k, l, a\}\subset kl$, then there are the first points $a_1\in[a, k]$ and $a_2\in [a, l]$ of the passes $g|_{[a, k]}$ and $g|_{[a, l]}$ respectively contained in $\partial U$.
	Since $g(kl) \cap g[i, p] \subset g(kl) \cap g(ij) = \varnothing$, we have $ga_1, ga_2 \in g[i, q] \subset g(ik)$.
	Denote by $b_1, b_2\in ik$ points such that $ga_1=gb_1$ and $ga_2=gb_2$.
	
	In this paragraph we show that arcs $I=[a_1, a_2]$, $\bar{I}=[b_1, b_2]$ satisfy the conditions of Lemma~\ref{l:LessLemma3}.
	Since $I\subset kl$, $\bar{I}\subset ki$, we have that $I, \bar{I}$ lie in $K_{3, 2}$.
	Since $g|_{ki}$ and $g|_{kl}$ are embeddings, we have that arcs $I$, $\bar{I}$ satisfy the property (\ref{l:LessLemma3.0}). 
	It follows from the definition of $a_1,a_2,b_1,b_2$ that the arcs $I$ and $\bar{I}$ satisfy property~(\ref{l:LessLemma3.1}).
	Let us show that arcs $I$, $\bar{I}$ satisfy the property (\ref{l:LessLemma3.2}).
	It follows from the definition of $a_1, a_2$ that $gI\cap U=\varnothing$.
	Since $gI\cap U=\varnothing$ and $g\bar{I}\subset \partial U$, we have $(gI\cup g\bar{I})\cap U=\varnothing$.
	Then there is an open component $W$ in $\R^2\backslash (gI\cup g\bar{I})$, such that $\Cl W\cap \Cl U = \Cl U$.
	Then $gk, gb, gl, gi \in \Cl W$, where $b$ is a point from $[5]\backslash\{i, j, k, l\}$. 
	Since $g(jb)\cap \partial W\subset g(jb)\cap (g(ki)\cup g(kl))=\varnothing$ and $gb\in \Cl W$, we have $gj\in \Cl W$.
	Then arcs $I$, $\bar{I}$ satisfy the property (\ref{l:LessLemma3.2}).
	
	Hence arcs $I, \bar{I}$ satisfy the conditions of Lemma \ref{l:LessLemma3}.
\end{proof} 
 
%\textbf{Proof of Lemma \ref{p:GreatLemma2} using  Lemmas \ref{l:improv_exist}, \ref{l:improv_exist_last_case}.} 

\section{Proof of Lemma \ref{p:GreatLemma2}}\label{s:proof}

It suffices to prove Lemma \ref{p:GreatLemma2} under the additional assumptions that 

$\bullet$ $g$ is a PL map in general position, and

%Figure \ref{ris:Reid} shows that for any PL almost embedding $g\colon \K \rightarrow \R^2$ there is an PL almost embedding $f\colon \K \rightarrow \R^2$ such that the restriction of $f$ to any edge of $K_{3,2}$ is an embedding and  $l(f) = l(g)$. 
%Then it suffices to prove Lemma \ref{p:GreatLemma2} under further additional assumption that restriction of $g$ to any edge of $K_{3,2}$ is embedding.

$\bullet$  $g|_{K_{3,2}}$ is not a PL embedding.

It suffices to prove that if $g|_{K_{3,2}}$ is not a PL embedding, then there is an improvement of $g$.  

\textit{Case 0:} {\it there is an edge $\alpha$ in $K_{3,2}$ such that $g|_{\alpha}$ is not an embedding.}
Figure \ref{ris:Reid} shows how to improve $g$.
%Since $\alpha$ is not in $123$, we have $l(f)=l(g)$. 

\textit{Case 1:} {\it the restriction of $g$ to any edge of $K_{3,2}$ is an embedding and there are no self-intersections of $g|_{ij\cup ik}$ for any $i \in \{4, 5\}$, $j\neq k\in [3]$.}
Since $g|_{K_{3,2}}$ is not a PL embedding, and the restriction of $g$ to any edge of $K_{3,2}$ is an embedding, and there are no self-intersections of $g|_{ij\cup ik}$ for any $i \in \{4, 5\}$, $j\neq k\in [3]$, there is a self-intersection of $g|_{i4\cup i5}$ for some $i \in [3]$.
Without loss of generality, assume that $g(14)\cap g(15)\varsupsetneq \{g1\}$.
Then  $(1, 5, 4)$ is interesting triple.
Denote by $p\in 15$ and $q\in 14$ the improving points for the interesting triple $(1, 5, 4)$.
By Lemma~\ref{l:improv} for the interesting triple $(1, 5, 4)$ there is a unique open component $U$ improving for the interesting triple $(1, 5, 4)$.
Moreover, we have $g1\in \partial U$.
In the following paragraph we prove that $g[5] \subset \Cl U$.

If $\{g4, g5\}\not\subset U$, then for some $i\in\{4, 5\}$ we have $\varnothing \neq g(i2)\cap \partial U \subset g(i2)\cap g(i1)$.
This contradicts the assumption that there are no self-intersections of $g|_{ij\cup ik}$ for any $i \in \{4, 5\}$, and $j\neq k\in [3]$.
Hence $g[5]\subset \Cl U$.

It follows that intervals $I:=[1, p], \bar{I}:=[1, q]$ satisfy the property (\ref{l:LessLemma3.2}).
By Lemma~\ref{l:improv} for the interesting triple $(1, 5, 4)$, we have that arcs $I$ and $\bar{I}$ lie in $K_{3, 2}$ and satisfy the properties (\ref{l:LessLemma3.0}), (\ref{l:LessLemma3.1}).
So we are done by Lemma \ref{l:LessLemma3}.

\textit{Case 2:} {\it the restriction of $g$ to any edge of $K_{3,2}$ is an embedding and there is a self-intersection of $g|_{ij\cup ik}$ for some $i \in \{4, 5\}$, $j\neq k\in [3]$.}
Without loss of generality, assume that $g(51) \cap g(52)\varsupsetneq \{g5\}$.
Then  $(5, 1, 2)$ and $(5, 2, 1)$ are interesting triples.
By Lemma~\ref{l:improv} for the interesting triple $(5, 1, 2)$ there is a unique open component $U$ improving for the interesting triple $(5, 1, 2)$.
%Moreover, we have $g5\in \partial U$.
By Lemma~\ref{l:improv} for the interesting triple $(5, 2, 1)$ there is a unique open component $U'$ improving for the interesting triple $(5, 2, 1)$.
%Moreover, we have $g5\in \partial U'$.
If $g1, g2\in U$, or $g2\notin U$, then from Lemma~\ref{l:improv_exist} for $i=5, j = 1$ and $ k = 2$ it follows that there are arcs satisfying the conditions of Lemma~\ref{l:LessLemma3}.
So we are done by Lemma \ref{l:LessLemma3}.
Hence it suffices to prove Lemma~\ref{p:GreatLemma2} under the additional assumption that $g2\in U$ and $g1\notin U$.
If $g2, g1\in U'$, or $g1\notin U'$, then from Lemma~\ref{l:improv_exist} for $i=5, j = 2$ and $ k = 1$ it follows that there are arcs satisfying the conditions of Lemma~\ref{l:LessLemma3}.
So we are done by Lemma \ref{l:LessLemma3}.
Hence it suffices to prove Lemma~\ref{p:GreatLemma2} under the additional assumption that $g1\in U'$ and $g2\notin U'$.

Since $g2\notin U'$ and $g4\in U'$, we have $\varnothing\neq g(24)\cap \partial U'= g(24)\cap g[5, q]$.
Then $g(24)\cap g(25)\varsupsetneq \{g2\}$.
Hence $(2, 5, 4)$ is interesting triple. 
By Lemma~\ref{l:improv} for the interesting triple $(2, 5, 4)$ there is a unique open component $W$ improving for the interesting triple $(2, 5, 4)$.
If $g5, g4\in W$, or $g4\notin W$, then from Lemma~\ref{l:improv_exist} for $i=2, j = 5$ and $ k = 4$ it follows that there are arcs satisfying the conditions of Lemma~\ref{l:LessLemma3}.
So we are done by Lemma \ref{l:LessLemma3}. 
Then it suffices to prove Lemma~\ref{p:GreatLemma2} under the additional assumption that $g4\in W$ and $g5\notin W$. 

Since $g1\notin U$, $g2\in U$, $g5\notin W$ and $g4\in W$, then from Lemma~\ref{l:improv_exist_last_case} for $i=5, j = 1, k = 2$ and $l = 4$ it follows that there are arcs satisfying the conditions of Lemma~\ref{l:LessLemma3}.
So we are done by Lemma \ref{l:LessLemma3}. $\Box$

\end{document}